\documentclass[reqno,12pt]{amsart}

\usepackage{mathtools}
\mathtoolsset{showonlyrefs}
\usepackage{etoolbox}
\usepackage[utf8]{inputenc}
\usepackage[T1]{fontenc}
\usepackage{amsmath,amssymb,amsfonts,mathrsfs,amsbsy}
\usepackage{relsize}
\usepackage{enumitem}
\usepackage{anyfontsize}
\usepackage{tikz,xcolor}
\usepackage{bm}
\usepackage{comment}
\usepackage[margin=1in]{geometry}

\definecolor{lava}{rgb}{0.81,0.06,0.13}
\definecolor{Cblue}{rgb}{0.50,0.85,0.85}
\definecolor{lime}{HTML}{A6CE39}

\usepackage[colorlinks,linkcolor=lava,citecolor=blue,urlcolor=Cblue,hypertexnames=false]{hyperref}

\everymath{\displaystyle}

\theoremstyle{plain}
\newtheorem{theorem}{Theorem}[section]
\newtheorem{proposition}{Proposition}[section]
\newtheorem{lemma}{Lemma}[section]
\newtheorem{corollary}{Corollary}[section]

\theoremstyle{definition}
\newtheorem{definition}{Definition}[section]
\newtheorem{remark}{Remark}[section]

\numberwithin{equation}{section}

\newcommand{\R}{\mathbb{R}}
\newcommand{\Sg}{S_\gamma}

\newcommand{\gN}{\gamma_N}

\DeclareRobustCommand{\orcidicon}{%
  \begin{tikzpicture}
    \draw[lime,fill=lime] (0,0) circle [radius=0.16]
      node[white] {{\fontfamily{qag}\selectfont\tiny ID}};
    \draw[white,fill=white] (-0.0625,0.095) circle [radius=0.007];
  \end{tikzpicture}\hspace{-2mm}}

\foreach \x in {A, B}{\expandafter\xdef\csname orcid\x\endcsname{\noexpand\href{https://orcid.org/\csname orcidauthor\x\endcsname}
            {\noexpand\orcidicon}}
}

\begin{document}
\title[Uniqueness for a concave Hardy parabolic equation]
{Uniqueness and comparison for a singular parabolic equation with a concave nonlinearity}
\author[M.~Loayza \& M.~Majdoub]{Miguel Loayza\orcidA{} and Mohamed Majdoub\orcidB{}}
\address[M. Loayza]{Departamento de Matem\'atica, Universidade Federal de Pernambuco - UFPE, Recife, PE, Brazil}
\email{miguel@dmat.ufpe.br}
\urladdr{https://orcid.org/0000-0002-3810-8928}
\address[M.~Majdoub]{Department of Mathematics, College of Science, Imam Abdulrahman Bin
Faisal University (IAU), P.O.~Box 1982, Dammam, Saudi Arabia}
\address[M.~Majdoub]{Basic and Applied Scientific Research Center, IAU, P.O.~Box 1982, 31441,
Dammam, Saudi Arabia}
\email{mmajdoub@iau.edu.sa}
\email{med.majdoub@gmail.com} \email{mohamed.majdoub@fst.rnu.tn}
\urladdr{https://orcid.org/0000-0001-6038-1069}
\subjclass[2020]{Primary 35A02, 35K58; Secondary 35K67, 35B51, 35B40, 35C06, 35K08.}

\keywords{Hardy parabolic equation; concave nonlinearity; non-Lipschitz nonlinearity;
uniqueness;
comparison principle; non-uniqueness for zero initial data; self-similar solutions;
large-time behavior; heat kernel estimates.}

\date{\today}

\begin{abstract}
We consider the Cauchy problem for $\partial_t u-\Delta u=|x|^{-\gamma}u^q$ in $\R^N$,
with $0<q<1$, $0<\gamma<\min\{2,N\}$ and nonnegative initial data
$u_0\in C_0(\R^N)$. Since $s\mapsto s^q$ is not Lipschitz continuous at the origin, the
uniqueness of nonnegative solutions is a delicate issue. {For $u_0\not\equiv0$, it
was known for small values of $\gamma$ and, more recently, under the additional
condition $\gamma<N(1-q)$.} We prove that, for $u_0\not\equiv0$,
uniqueness holds on the whole range $0<\gamma<\min\{2,N\}$, which settles a conjecture
of the authors. In fact, we establish a comparison principle between bounded
subsolutions and supersolutions which are bounded from below by a multiple of the
profile $t^{1/(1-q)}(|x|+\sqrt t\,)^{-\gamma/(1-q)}$. We then describe completely the
set of nonnegative solutions emanating from $u_0\equiv0$: it consists of the trivial
solution, of a unique positive solution $U$, and of the delayed copies
$U(\cdot,(t-\tau)_+)$, $\tau>0$. The solution $U$ is radially symmetric and
self-similar, it obeys {two-sided pointwise bounds by multiples of the same
profile}, and it is the limit of the
solutions with nontrivial data as the data decrease to zero. Finally, we show that $U$
describes the large-time behavior of every nontrivial solution. The proofs replace
the singular Gronwall inequality used in earlier work by a quotient estimate which
exploits the concavity of $s\mapsto s^q$ and requires no smallness condition.
\end{abstract}

\maketitle

\section{Introduction and main results}\label{sec:intro}

We study the Cauchy problem for the semilinear heat equation with a singular weight and
a concave reaction term,
\begin{equation}\label{eq:main}
\begin{cases}
\partial_t u-\Delta u=|x|^{-\gamma}u^q & \text{in } \R^N\times(0,\infty),\\[1mm]
u(x,0)=u_0(x) & \text{in }\R^N,
\end{cases}
\qquad 0<q<1,\ \ \gamma>0,
\end{equation}
where $u_0\in C_0(\R^N)$ and $u_0\ge0$. Here $C_0(\R^N)$ denotes the closure of
$C_c^\infty(\R^N)$ in $L^\infty(\R^N)$, that is, the space of continuous functions
vanishing at infinity. The weight $|x|^{-\gamma}$ is of Hardy type, and the reaction
term $u^q$ is concave and sublinear.

Two thresholds restrict the admissible values of $\gamma$. First, the weight must be
locally integrable, which requires $\gamma<N$. Second, the value $\gamma=2$ is critical
for the heat operator. Indeed, for $\gamma=2$ the existence of solutions to
\eqref{eq:main} fails in general; see \cite[Remark 1.2(b)]{LM} and the references
therein, notably \cite{BG,CRB,HiTs}. The decisive role of the inverse-square potential
for linear heat equations goes back to Baras and Goldstein \cite{BG}, and it was further
{analyzed} in \cite{CM,VZ}. The natural range of parameters is therefore
\begin{equation}\label{eq:gN}
0<\gamma<\gN:=\min\{2,N\}.
\end{equation}

As usual, we work with the integral formulation of \eqref{eq:main}. Let
$\{S(t)\}_{t\ge0}$ denote the heat semigroup on $\R^N$, given by $S(t)\varphi=G_t\ast\varphi$
with $G_t(x)=(4\pi t)^{-N/2}e^{-|x|^2/4t}$ for $t>0$. We set
\begin{equation}\label{eq:duhamel}
u(t)=S(t)u_0+\int_0^t \Sg(t-\sigma)\,u^q(\sigma)\,d\sigma,
\qquad \Sg(t)\phi:=S(t)\big(|\cdot|^{-\gamma}\phi\big).
\end{equation}

\begin{definition}[\cite{LM}]\label{def:sol}
A \emph{nonnegative solution} of \eqref{eq:main} is a function
$u\in C([0,\infty);C_0(\R^N))$ with $u\ge0$ which satisfies \eqref{eq:duhamel} for every
$t\ge0$.
\end{definition}

For $0<\gamma<\gN$ and every $0\le u_0\in C_0(\R^N)$, problem \eqref{eq:main} admits a
nonnegative solution in the sense of Definition~\ref{def:sol}; see
\cite[Theorem~1.1]{LM}. The present paper is concerned with two questions which were
left open in \cite{LM}: the uniqueness of this solution when $u_0\not\equiv0$, and the
description of all solutions when $u_0\equiv0$.

\subsection*{Background}
For superlinear powers, the equation $\partial_t u-\Delta u=|x|^{-\gamma}|u|^{p-1}u$ with
$p>1$ is known as the Hardy parabolic equation. It belongs to the wider family of
Hardy--H\'enon parabolic equations, which has been studied intensively in recent years.
Well-posedness in Lebesgue spaces was established by Ben Slimene, Tayachi and Weissler
\cite{BSTW}, extending the classical theory of Weissler \cite{W80} for $\gamma=0$.
Uniqueness and non-uniqueness in critical spaces were investigated by Tayachi
\cite{Tayachi} and by Chikami, Ikeda, Taniguchi and Tayachi \cite{CITT}. The optimal
singularities of initial data for local solvability were identified by Hisa and Takahashi
\cite{HT} and by Hisa and Sier\.z\c{e}ga \cite{HS}. Local existence for singular or uniformly
local initial data was obtained in \cite{CGL,CCL}, where more general nonlinearities are
also allowed. Equations with inhomogeneous nonlinearities, forcing terms or stochastic
perturbations are treated, for instance, in \cite{Maj21,Maj23,MM21}. We refer to the
monograph of Quittner and Souplet \cite{QS} for a comprehensive account of superlinear
parabolic problems.

The sublinear case $0<q<1$ has a different nature. For $\gamma=0$, the Cauchy problem
was studied in detail by Aguirre and Escobedo \cite{AE}. On bounded domains, equations
of the form $\partial_t u-\Delta u=a(x)u^q$ were considered in the pioneering work of
Fujita and Watanabe \cite{FW} for continuous $a$, and by Loayza \cite{Loayza06} for
singular $a$ and singular initial data. Concave--convex nonlinearities
$\lambda u^q+u^p$ on bounded domains were {analyzed} by Cazenave, Dickstein and Escobedo
\cite{CDE}. In the appendix of \cite{CDE}, the existence and uniqueness of positive
solutions are established for a class of non-Lipschitz nonlinearities, following
the ideas of \cite{FW} and of Escobedo and Herrero \cite{EH93}.

\subsection*{The uniqueness problem}
The main difficulty is the lack of uniqueness inherent in sublinear reaction terms.
The map $s\mapsto s^q$ is not Lipschitz continuous at $s=0$. Moreover,
$\int_{0^+}s^{-q}\,ds<\infty$, so that Osgood-type uniqueness criteria
\cite{AgarwalLakshmikantham,Osgood,Nagumo} cannot be applied near the set $\{u=0\}$.
This is not a technical {artifact}. For $\gamma=0$ and $u_0\equiv0$, the problem possesses,
besides the trivial solution, the one-parameter family of solutions
\[
u_s(x,t)=\big[(1-q)(t-s)_+\big]^{1/(1-q)},\qquad s\ge0,
\]
see \cite{AE}; an analogous family exists for reaction--diffusion systems with
non-Lipschitz coupling \cite{EH}. On bounded domains, every solution with zero initial
data is a time translate of the unique positive one \cite[Remark~6.4]{CDE}. Non-uniqueness phenomena of a different nature, caused
by singular initial data and superlinear growth, are also classical; see Haraux and
Weissler \cite{HW} and Br\'ezis and Cazenave \cite{BC}.

The assumption $u_0\not\equiv0$ prevents the solution from staying on the degenerate set
$\{u=0\}$. For $\gamma=0$, Aguirre and Escobedo \cite{AE} resolved the uniqueness
issue, and a key ingredient of their analysis is the fact that every positive solution
satisfies $u(t)\ge[(1-q)t]^{1/(1-q)}$; see also \cite[(1-2) and Remark~1.4]{LM}. For $\gamma>0$, the authors proved in \cite[Theorem~1.3(i)]{LM} that every
nonnegative solution with initial data $0\not\equiv u_0\ge0$ satisfies the pointwise
lower bound
\begin{equation}\label{eq:lower}
u(x,t)\ \ge\ w(x,t):=\big[(1-q)\eta_0\big]^{1/(1-q)}\,t^{1/(1-q)}\,
\big(|x|+\sqrt t\,\big)^{-\theta},
\qquad \theta:=\frac{\gamma}{1-q},
\end{equation}
for all $x\in\R^N$ and $t\ge0$, where
\begin{equation}\label{eq:eta0}
\eta_0=\eta_0(\gamma)=(4\pi)^{-N/2}\int_{\R^N}e^{-|z|^2/4}\,(1+|z|)^{-\gamma/(1-q)}\,dz
\end{equation}
is the constant defined in \cite[(1-6)]{LM}. The exponent $\theta=\gamma/(1-q)$ is dictated
by the structure of the equation. Indeed, if $u$ behaved like $|x|^{-\theta}$, the
reaction term would behave like $|x|^{-\gamma-\theta q}$. Balancing the two exponents gives
$\gamma+\theta q=\theta$, and we obtain $\theta=\gamma/(1-q)$. Throughout the paper we
write
\begin{equation}\label{eq:W}
W(x,t):=t^{1/(1-q)}\big(|x|+\sqrt t\,\big)^{-\theta},\qquad
c_0:=\big[(1-q)\eta_0\big]^{1/(1-q)},\qquad
\alpha:=\frac{2-\gamma}{1-q}.
\end{equation}
The profile $W$ is invariant under the natural scaling of the
equation, namely $W(\lambda x,\lambda^2t)=\lambda^{\alpha}W(x,t)$ for $\lambda>0$, and
$\sup_{x\in\R^N}W(x,t)=W(0,t)=t^{\alpha/2}$.

Using \eqref{eq:lower}, it was shown in \cite[Theorem~1.3(ii)]{LM} that uniqueness
holds for $0<\gamma<\gamma^*$, for some $\gamma^*\in(0,\gN)$. The argument compares two
solutions $u$ and $v$ through the positive part $z=[v-u]_+$ of their difference. By the
mean value theorem and \eqref{eq:lower}, this function satisfies
\[
z(x,t)\le\frac{q}{(1-q)\eta_0}\int_0^t\sigma^{-1}
\Big[\Sg(t-\sigma)\big((|\cdot|+\sqrt\sigma)^\gamma z(\sigma)\big)\Big](x)\,d\sigma,
\]
and an iteration based on a singular Gronwall inequality \cite{Henry} yields the decay
estimate $\|z(t)\|_\infty\le\Lambda(\gamma)^k\,C(t)$ for every $k\ge1$. The factor
$\Lambda(\gamma)$ tends to $q<1$ as $\gamma\to0$, which gives uniqueness for small
$\gamma$. However, $\Lambda(\gamma)$ blows up as $\gamma\uparrow\gN$. This blow-up has two
sources. The temporal factor $(t-\sigma)^{-\gamma/2}$ produced by the smoothing of the
weight ceases to be integrable as $\gamma\to2$, and the weight $|x|^{-\gamma}$ ceases to
be locally integrable as $\gamma\to N$. The authors conjectured \cite[Remark~4.3]{LM}
that uniqueness holds on the whole range $0\le\gamma<\gN$. {Recently, Ioku and
Nitta \cite{IokuNitta2026} proved uniqueness for $0<\gamma<\min\{\gN,N(1-q)\}$ by an
argument based on the ratio of two solutions; we compare their work with ours in
Remark~\ref{rem:IokuNitta} below. Their result confirms the conjecture when $N\ge3$
and $q\le1-2/N$, and leaves open the range $N(1-q)\le\gamma<\gN$, which is nonempty
when $N\in\{1,2\}$, and when $N\ge3$ and $q>1-2/N$.}

\subsection*{Solutions emanating from zero}
When uniqueness fails, the natural question is to describe the whole set of solutions.
For non-Lipschitz reaction terms this has been achieved in three settings. For
$\gamma=0$ in $\R^N$, Aguirre and Escobedo \cite[Corollary~2.6]{AE} proved that the
nontrivial nonnegative mild solutions with $u_0\equiv0$ are exactly the functions
$u_s$ displayed above; in other words, they are the positive solution
$[(1-q)t]^{1/(1-q)}$ and its delayed copies. For the weakly coupled system
$u_t-\Delta u=v^p$, $v_t-\Delta v=u^q$ with $pq<1$ in a bounded domain with Dirichlet
boundary conditions, Escobedo and Herrero \cite[Theorem~1(b)]{EH93} showed that the
solutions emanating from $(0,0)$ are the trivial one, a solution $(U,V)$ which is
positive for $t>0$, and the one-parameter family $(U,V)(\cdot,(t-\mu)_+)$, $\mu>0$.
Dickstein and Escobedo \cite[Theorem~1.1]{DE} extended this trichotomy to the strongly
coupled system with product nonlinearities $u^{a}v^{b}$ and $u^{c}v^{d}$ in the
non-Lipschitz regime. In all three cases the structure is the same: the solution
may rest at zero for an arbitrary time, and once it leaves zero its evolution is
uniquely determined.

For the weighted problem \eqref{eq:main} with $\gamma>0$ this question was not
addressed in \cite{LM}. Two features distinguish it from the situations above. In
contrast with the case $\gamma=0$, no explicit positive solution is available, and the
equation is not invariant under spatial translations, so that the solutions leaving
zero cannot be independent of $x$. In contrast with bounded domains, the lower bound
\eqref{eq:lower} degenerates as $|x|\to\infty$, and there is no first eigenfunction to
compare with. We show below that the same trichotomy nevertheless holds, and that the
positive solution leaving zero can be described rather precisely.

\subsection*{Main results}
Our approach avoids the {linearization} of the difference of two solutions, and hence
the singular Gronwall inequality. Instead, it relies on a comparison argument which
exploits the concavity of $s\mapsto s^q$ in a global way. Our first result is the
following comparison principle, which proves the conjecture of \cite{LM}.

\begin{theorem}\label{thm:comparison}
Let $N\ge1$, $0<q<1$ and $0<\gamma<\gN$. Let $u_0,v_0\in C_0(\R^N)$ satisfy
$0\le u_0\le v_0$ and $v_0\not\equiv0$, and let $v$ be a nonnegative solution of
\eqref{eq:main} with initial data $v_0$. Let $u:\R^N\times(0,\infty)\to[0,\infty)$ be a
measurable function which is bounded on $\R^N\times(0,T]$ for every $T>0$ and satisfies
\begin{equation}\label{eq:subsol}
u(t)\le S(t)u_0+\int_0^t\Sg(t-\sigma)\,u^q(\sigma)\,d\sigma
\qquad\text{in }\R^N,\ \text{for every }t>0.
\end{equation}
Then $u\le v$ in $\R^N\times(0,\infty)$.
\end{theorem}

In Theorem~\ref{thm:comparison} the initial data $u_0$ may vanish identically, and $u$
need not be continuous. When $u$ is itself a solution, we obtain the following
consequence.

\begin{corollary}\label{cor:uniq}
Let $N\ge1$, $0<q<1$ and $0<\gamma<\gN$.
\begin{enumerate}[label=\rm(\roman*)]
\item If $0\le u_0\in C_0(\R^N)$ and $u_0\not\equiv0$, then problem \eqref{eq:main}
has exactly one nonnegative solution.
\item If $u$ and $v$ are nonnegative solutions of \eqref{eq:main} with initial data
$0\le u_0\le v_0$ and $v_0\not\equiv0$, then $u\le v$ in $\R^N\times[0,\infty)$.
\end{enumerate}
\end{corollary}
\begin{remark}\rm \label{rem:IokuNitta}
Ioku and Nitta \cite{IokuNitta2026} recently considered the equation
$\partial_tu-\Delta u=a(x)|x|^{-\gamma}|u|^{q-1}u$ with $a=\mathbf 1_{B_L}$, where
$B_L$ is the open ball of radius $0<L\le+\infty$ centered at the origin. For $L=+\infty$
this is \eqref{eq:main}. They proved \cite[Theorem~1.1]{IokuNitta2026} that, for
$0\le u_0\in C_0(\R^N)$ with $u_0\not\equiv0$, the problem has a unique solution in
$C([0,\infty);C_0(\R^N))$, provided that $0<\gamma<\min\{\gN,N(1-q)\}$. Their argument,
like ours, replaces the difference of two solutions by their ratio. If $u$ and $v$
are two solutions with the same initial data and $K^{-1}v\le u\le Kv$ on
$B_L\times(0,t]$ for some $K\ge1$, the homogeneity of $s\mapsto s^q$ gives the same
inequalities with $K$ replaced by $K^q$, and therefore $K=1$. The main step is to show
that such a constant exists near $t=0$. To this end, Ioku and Nitta combine a lower
bound of the type \eqref{eq:lower} \cite[Lemma~2.2]{IokuNitta2026}, which for
$L=+\infty$ is \eqref{eq:lower} itself, taken from \cite{LM}, with an upper bound of the form
$u\le\varphi_\infty+C\,t^{1/(1-q)}|x|^{-\gamma/(1-q)}$, where $\varphi_\infty$ is the
minimal solution \cite[Lemma~2.4]{IokuNitta2026}. This upper bound is obtained by
iterating the integral equation and requires the local integrability of
$|x|^{-\gamma/(1-q)}$. This is the origin of the condition $\gamma<N(1-q)$
\cite[Remark~2.5]{IokuNitta2026}.
\end{remark}
\begin{remark}\rm 
The present paper differs from \cite{IokuNitta2026} in four respects. 
\begin{enumerate}[label=\rm(\roman*)]
\item First,
Corollary~\ref{cor:uniq}(i) holds on the whole range $0<\gamma<\gN$. The condition
$\gamma<N(1-q)$ is a genuine restriction when $N\in\{1,2\}$, for every $q\in(0,1)$, and
when $N\ge3$ and $q>1-2/N$. In these cases, uniqueness for $N(1-q)\le\gamma<\gN$ is new.
When $N\ge3$ and $q\le1-2/N$, the uniqueness statement of Corollary~\ref{cor:uniq}(i)
is contained in \cite{IokuNitta2026}, and the contribution of the present paper lies in
the three points below.
\item Second, we prove comparison principles and not only uniqueness. In
Theorem~\ref{thm:comparison}, the function $u$ is an arbitrary bounded measurable
subsolution whose initial data lies below that of $v$, and in Theorem~\ref{thm:super}
the solution $v$ is replaced by any supersolution lying above a multiple of $W$. The
quotient in Lemma~\ref{lemma.Comp} is one-sided: it controls $(u-v)_+/v$ only, and it
requires neither an upper bound for $v$ nor a lower bound for $u$. Its finiteness is
obtained without any upper estimate for solutions, by first comparing $u$ with the
supersolution $Bv+\varepsilon$. This step uses only the bound of Lemma~\ref{lemma.H}
on the heat flow of $|x|^{-\gamma}$, which is finite as soon as $\gamma<N$. This is why
no condition other than $\gamma<\gN$ appears in our results.
\item Third, we describe all the solutions emanating from $u_0\equiv0$
(Theorem~\ref{thm:zero}), a case not considered in \cite{IokuNitta2026}. It cannot be
handled by a two-sided ratio estimate between solutions with the same initial datum,
since $U$ and $U_\tau$, $\tau>0$, have the same initial data and $U/U_\tau=\infty$ on
$\R^N\times(0,\tau]$. Our argument rests instead on Theorem~\ref{thm:super}. It
produces a self-similar solution $U$ with explicit two-sided bounds, which is the
largest bounded subsolution with zero initial data and the limit of the solutions with
vanishing nontrivial data. The supersolution $K_0W$ provided by Lemma~\ref{lemma.W} is
a new ingredient of this part.
\item Fourth, Theorem~\ref{thm:asym}(i) gives explicit two-sided bounds for every
nontrivial solution, whereas the upper bounds of \cite[Lemma~2.4]{IokuNitta2026} are
expressed in terms of the minimal solution $\varphi_\infty$. The large-time behavior
described in Theorem~\ref{thm:asym}(ii) has no counterpart in \cite{IokuNitta2026} or
\cite{LM}.
\end{enumerate}
\end{remark}
In the proof of Theorem~\ref{thm:comparison}, the assumption $v_0\not\equiv0$ is used
only to guarantee the lower bound \eqref{eq:lower} for $v$. The argument therefore
yields a more flexible statement, in which $v$ is merely a supersolution lying above
a multiple of the profile $W$. This version is the tool which allows us to treat the
zero initial datum.

\begin{theorem}\label{thm:super}
Let $N\ge1$, $0<q<1$, $0<\gamma<\gN$, $T>0$ and $0\le u_0\in C_0(\R^N)$. Let
$u:\R^N\times(0,T]\to[0,\infty)$ and $v:\R^N\times(0,T]\to(0,\infty)$ be measurable
functions such that, for every $t\in(0,T]$,
\begin{equation}\label{eq:subsuper}
u(t)\le S(t)u_0+\int_0^t \Sg(t-s)\, u^q(s)\, ds
\quad\text{and}\quad
v(t)\geq S(t)u_0+\int_0^t \Sg(t-s)\, v^q(s)\, ds
\quad\text{in }\R^N.
\end{equation}
Assume that $u$ is bounded and that there exists a constant $c>0$ such that
\begin{equation}\label{eq:vW}
v\ge c\,W\qquad{\text{in }\R^N\times(0,T],}
\end{equation}
where $W$ is given by~\eqref{eq:W}. 
Then $u\le v$ in $\R^N\times(0,T]$.
\end{theorem}

We now turn to the zero initial datum. For $\tau\ge0$ and a function $U$ defined on
$\R^N\times[0,\infty)$, we write $U_\tau(x,t):=U(x,(t-\tau)_+)$.

\begin{theorem}\label{thm:zero}
Let $N\ge1$, $0<q<1$ and $0<\gamma<\gN$, and let $W$, $c_0$ and $\alpha$ be given by
\eqref{eq:W}.
\begin{enumerate}[label=\rm(\roman*)]
\item Problem \eqref{eq:main} with $u_0\equiv0$ has a unique nonnegative solution $U$
such that $U(t)\not\equiv0$ for every $t>0$. This solution is radially symmetric with
respect to $x$ and self-similar,
\begin{equation}\label{eq:selfsim}
U(x,t)=t^{\alpha/2}\,\Phi\Big(\frac{x}{\sqrt t}\Big),\qquad \Phi:=U(\cdot,1)\in C_0(\R^N),
\end{equation}
and there exists a constant $K_0=K_0(N,\gamma,q)>0$ such that
\begin{equation}\label{eq:Ubounds}
c_0\,W\le U\le K_0\,W\qquad{\text{in }\R^N\times[0,\infty).}
\end{equation}
Equivalently, $c_0(1+|z|)^{-\theta}\le\Phi(z)\le K_0(1+|z|)^{-\theta}$ for all
$z\in\R^N$.
\item The nonnegative solutions of \eqref{eq:main} with $u_0\equiv0$ are exactly the
trivial solution $u\equiv0$ and the functions $U_\tau$, $\tau\ge0$.
\item The solution $U$ is the largest element of the set of nonnegative measurable
functions on $\R^N\times(0,\infty)$ which are bounded on $\R^N\times(0,T]$ for every
$T>0$ and satisfy \eqref{eq:subsol} with $u_0\equiv0$. Moreover, $U\le v$ for every
nonnegative solution $v$ of \eqref{eq:main} with initial data $v_0\not\equiv0$.
\item Let $(v_{0,n})_{n\ge1}\subset C_0(\R^N)$ be a nonincreasing sequence of nonnegative
functions, none of which vanishes identically, such that $v_{0,n}\to0$ pointwise, and
let $v_n$ be the solution of \eqref{eq:main} with initial data $v_{0,n}$. Then
$v_n\downarrow U$ uniformly on $\R^N\times[0,T]$ for every $T>0$.
\end{enumerate}
\end{theorem}

Theorem~\ref{thm:zero}(ii) is the counterpart for \eqref{eq:main} of
\cite[Corollary~2.6]{AE}, \cite[Theorem~1(b)]{EH93} and \cite[Theorem~1.1(iii)]{DE}.
{Parts (i), (iii) and (iv) describe features which are absent for $\gamma=0$: the
weight forces the positive solution leaving zero to decay like $|x|^{-\theta}$ at
infinity and to be self-similar with a nonconstant profile, whereas the positive
solution of Aguirre and Escobedo is constant in space. For $\gamma=0$, there is no
nontrivial solution in $C_0(\R^N)$ at all; see Remark~\ref{rem:main}(b).}

Our last result shows that $U$ governs the dynamics of all nontrivial solutions.

\begin{theorem}\label{thm:asym}
Let $N\ge1$, $0<q<1$ and $0<\gamma<\gN$. Let $0\le u_0\in C_0(\R^N)$ with
$u_0\not\equiv0$, and let $u$ be the solution of \eqref{eq:main} with initial datum
$u_0$.
\begin{enumerate}[label=\rm(\roman*)]
\item There exists a constant $K_1=K_1(N,\gamma,q)>0$ such that
\begin{equation}\label{eq:twosided}
U\le u\le 2\|u_0\|_\infty+K_1\,W\qquad\text{in }\R^N\times[0,\infty).
\end{equation}
In particular, $c_0\,t^{\alpha/2}\le\|u(t)\|_\infty\le2\|u_0\|_\infty+K_1\,t^{\alpha/2}$
for every $t\ge0$.
\item We have
\begin{equation}\label{eq:asym}
\lim_{t\to\infty}t^{-\alpha/2}\,\big\|u(t)-U(t)\big\|_\infty=0,
\end{equation}
that is, $t^{-\alpha/2}u(\sqrt t\,z,t)\to\Phi(z)$ uniformly with respect to $z\in\R^N$ as
$t\to\infty$.
\end{enumerate}
\end{theorem}

\begin{remark}\rm \label{rem:main}
~\begin{enumerate}[label=\rm(\alph*)]
\item The proofs of Theorems~\ref{thm:comparison} and~\ref{thm:super} use the structure
of the problem only through three facts: the concavity of $s\mapsto s^q$, the lower
bound \eqref{eq:lower} or its substitute \eqref{eq:vW}, and an upper bound for the heat
flow of the weight (Lemma~\ref{lemma.H} below). The last two estimates involve the
same space-time profile $W$, and this matching is the heart of the argument. Both
thresholds in \eqref{eq:gN} are used: $\gamma<N$ makes the weight locally integrable,
and $\gamma<2$ makes its heat flow integrable in time.
\item {The case $\gamma=0$ requires a change of functional framework. Indeed,
the lower bound $u(t)\ge[(1-q)t]^{1/(1-q)}$ of \cite{AE} shows that for $\gamma=0$ no
nontrivial solution remains in $C_0(\R^N)$, and this is why \cite{AE} works in weighted
spaces of functions which may grow at infinity. On the other hand,
Lemma~\ref{lemma.Comp} does not involve the space $C_0(\R^N)$, and for $\gamma=0$ one
has $H(x,t)=t$ and $W(x,t)=t^{1/(1-q)}$. Hence the proof of Theorem~\ref{thm:super}
applies verbatim to bounded mild solutions with bounded nonnegative initial data, and
it recovers, in this class, the uniqueness result of \cite{AE} together with a
comparison principle. In this sense uniqueness holds on the full range
$0\le\gamma<\gN$ predicted in \cite[Remark~4.3]{LM}. The solutions of \cite{AE} with
unbounded data are not covered, since our argument uses the boundedness of the
subsolution.}
\item In the argument, the contraction is provided by the inequality $\Theta\le q\Theta$
of Lemma~\ref{lemma.Comp}, which holds for every $\gamma$. {Note that} $q$
is exactly the limit of the factor $\Lambda(\gamma)$ of \cite{LM} as $\gamma\to0$.
\item The condition $v_0\not\equiv0$ cannot be removed from
Corollary~\ref{cor:uniq}(ii). Indeed, by Theorem~\ref{thm:zero}(ii) the solutions $U$
and $U_\tau$ with $\tau>0$ have the same initial datum, and $U\not\le U_\tau$ because
$U_\tau(t)\equiv0<U(t)$ for $0<t\le\tau$. Note that the family
$\{0\}\cup\{U_\tau\}_{\tau\ge0}$ is totally ordered: by
Theorem~\ref{thm:zero}(iii) we have $U_\tau\le U$ for every $\tau\ge0$, {and if
$0\le\sigma\le\tau$, the inequality $U_{\tau-\sigma}\le U$ evaluated at the time
$t-\sigma$ gives $U_\tau(t)\le U_\sigma(t)$ for $t\ge\sigma$, while both sides vanish
for $t<\sigma$.}
\item {Inserting \eqref{eq:selfsim} into the equation shows, at least
formally, that the profile $\Phi$ is a positive radial solution of
\[
-\Delta\Phi-\frac{z}{2}\cdot\nabla\Phi+\frac{\alpha}{2}\,\Phi=|z|^{-\gamma}\,\Phi^q
\qquad\text{in }\R^N\setminus\{0\}.
\]
We do not study this elliptic problem here. We expect $\Phi$ to be radially
nonincreasing and to satisfy $|z|^{\theta}\Phi(z)\to\ell$ as $|z|\to\infty$ for some
$\ell\in[c_0,K_0]$.}
\item For $\gamma=0$, the large-time behavior of the solutions of \cite{AE} depends on
the behavior of $u_0$ at infinity and is described by a family of self-similar
solutions \cite[Section~3]{AE}. In the present setting, the initial data belong to
$C_0(\R^N)$ and Theorem~\ref{thm:asym}(ii) shows that all nontrivial solutions share
the same asymptotic profile $\Phi$. Since $\Phi\ge c_0(1+R)^{-\theta}$ on the ball
$\{|z|\le R\}$, the limit \eqref{eq:asym} also gives $u/U\to1$ as $t\to\infty$,
uniformly on the parabolic regions $\{|x|\le R\sqrt t\}$, $R>0$.
\end{enumerate}
\end{remark}

\subsection*{{Organization} and notation}
Section~\ref{sec:tools} collects the tools needed in the proofs: the invariances of
the integral equation, the lower bound \eqref{eq:lower}, two estimates for the heat
flow of the weight $|x|^{-\gamma}$, and a regularity property of the Duhamel integral.
Section~\ref{sec:comp} contains the comparison lemma on which the whole argument rests.
Theorems~\ref{thm:super} and~\ref{thm:comparison} and Corollary~\ref{cor:uniq} are
proved in Section~\ref{sec:proofs}. Section~\ref{sec:zero} is devoted to the solutions
emanating from zero and contains the proof of Theorem~\ref{thm:zero}, and
Theorem~\ref{thm:asym} is proved in Section~\ref{sec:asym}. Section~\ref{sec:concl}
contains some concluding remarks. Throughout the paper, $\omega_N$ denotes the volume
of the unit ball of $\R^N$, so that the surface measure of the unit sphere equals
$N\omega_N$. For $a\in\R$ we write $a_+=\max\{a,0\}$. All functions are real valued,
and inequalities between functions are understood pointwise. We denote by
$BUC(\R^N)$ the space of bounded uniformly continuous functions on $\R^N$, endowed
with the supremum norm.

\section{Preliminaries}\label{sec:tools}

\subsection{The heat semigroup and the Duhamel integrals}
The heat semigroup is positivity preserving and satisfies $S(t)1=1$ for every $t\ge0$.
Consequently, for every nonnegative measurable function $f$ on $\R^N\times(0,t)$, the
Duhamel integral
\[
\int_0^t\Sg(t-\sigma)f(\sigma)\,d\sigma\,(x)
=\int_0^t\int_{\R^N}G_{t-\sigma}(x-y)\,|y|^{-\gamma}f(y,\sigma)\,dy\,d\sigma
\]
is well defined as an element of $[0,\infty]$ by Tonelli's theorem. It is monotone and
positively homogeneous with respect to $f$. These are the only properties of $\Sg$ used
in Section~\ref{sec:comp}.

The integral equation \eqref{eq:duhamel} inherits the invariances of the differential
equation: translations in time, rotations, and the scaling associated with the
exponent $\alpha$ of \eqref{eq:W}.

\begin{lemma}\label{lemma.inv}
Let $0<q<1$ and $0<\gamma<\gN$, and let $u$ be a nonnegative solution of
\eqref{eq:main} with initial data $u_0$.
\begin{enumerate}[label=\rm(\roman*)]
\item For every $s\ge0$, the function $(x,t)\mapsto u(x,t+s)$ is a nonnegative solution
of \eqref{eq:main} with initial data $u(s)$.
\item For every $\lambda>0$, the function
$u_\lambda(x,t):=\lambda^{-\alpha}u(\lambda x,\lambda^2t)$ is a nonnegative solution of
\eqref{eq:main} with initial data $\lambda^{-\alpha}u_0(\lambda\,\cdot)$.
\item For every orthogonal matrix $R$, the function $(x,t)\mapsto u(Rx,t)$ is a
nonnegative solution of \eqref{eq:main} with initial data $u_0(R\,\cdot)$.
\end{enumerate}
\end{lemma}

\begin{proof}
In each case the new function clearly belongs to $C([0,\infty);C_0(\R^N))$, and we
only have to check \eqref{eq:duhamel}. Set $g(\sigma):=|\cdot|^{-\gamma}u^q(\sigma)\ge0$.
For (i), we split the Duhamel integral in \eqref{eq:duhamel} at time $t+s$ as
$\int_0^s+\int_s^{t+s}$. Since $G_t\ast G_{s-\sigma}=G_{t+s-\sigma}$, Tonelli's theorem
gives
\[
S(t+s)u_0+\int_0^sS(t+s-\sigma)g(\sigma)\,d\sigma
=S(t)\Big[S(s)u_0+\int_0^sS(s-\sigma)g(\sigma)\,d\sigma\Big]=S(t)u(s),
\]
while the change of variable $\sigma=s+\sigma'$ transforms the remaining integral into
$\int_0^t\Sg(t-\sigma')u^q(s+\sigma')\,d\sigma'$. This proves (i). For (ii), we
evaluate \eqref{eq:duhamel} at the point $(\lambda x,\lambda^2t)$ and multiply by
$\lambda^{-\alpha}$. The identity $G_{\lambda^2t}(\lambda z)=\lambda^{-N}G_t(z)$ gives
$[S(\lambda^2t)u_0](\lambda x)=[S(t)u_0(\lambda\,\cdot)](x)$. In the Duhamel integral
we substitute $y=\lambda y'$ and $\sigma=\lambda^2\sigma'$, and we use
$u^q(\lambda y',\lambda^2\sigma')=\lambda^{\alpha q}u_\lambda^q(y',\sigma')$. The
resulting power of $\lambda$ is $\lambda^{-\alpha-N+N+2-\gamma+\alpha q}=1$, because
$\alpha(1-q)=2-\gamma$. This proves (ii). Assertion (iii) follows from the change of
variable $y=Ry'$, since $G_t$ and the weight are radial.
\end{proof}

\subsection{A lower bound for nontrivial solutions}
We recall the lower bound established in \cite{LM}.

\begin{proposition}[{\cite[Theorem~1.3(i)]{LM}}]\label{prop:lower}
Let $0<q<1$, $0<\gamma<\gN$ and $0\le u_0\in C_0(\R^N)$ with $u_0\not\equiv0$. Then every
nonnegative solution $u$ of \eqref{eq:main} satisfies \eqref{eq:lower} for all
$x\in\R^N$ and $t\ge0$, with $\eta_0$ given by \eqref{eq:eta0}. In particular, $u>0$
in $\R^N\times(0,\infty)$.
\end{proposition}

\begin{remark}\rm \label{rem:local}
The proof of \cite[Theorem~1.3(i)]{LM} uses the integral equation
\eqref{eq:duhamel} only on the time interval $[0,t]$. It first treats initial data
bounded from below by a Gaussian, through an iteration of the inequality
$u(t)\ge\int_0^t\Sg(t-\sigma)u^q(\sigma)\,d\sigma$, and then reduces the general case to
this one by applying the result to $u(\cdot+\varepsilon)$ and letting
$\varepsilon\to0$. Therefore \eqref{eq:lower} also holds on $\R^N\times[0,T]$ for
solutions defined on a finite interval $[0,T]$, and Theorem~\ref{thm:comparison} and
Corollary~\ref{cor:uniq} extend accordingly to such solutions.
\end{remark}

\subsection{The heat flow of the weight}
The following estimate shows that the Duhamel integral of the constant function $1$ has
exactly the profile appearing in \eqref{eq:lower}, once it is raised to the power
$1/(1-q)$.

\begin{lemma}\label{lemma.H}
Let $0<\gamma<\gN$ and set
\begin{equation}\label{H}
H(t):=\int_0^t S(t-s)|\cdot|^{-\gamma}\, ds,\qquad t>0.
\end{equation}
There exists a constant $C_1=C_1(N,\gamma)>0$ such that
\[
H(x,t)\le C_1\, t\, \big(|x|+\sqrt{t}\,\big)^{-\gamma}
\qquad\text{for all } x\in \R^N \text{ and } t>0.
\]
\end{lemma}

\begin{proof}
Since $|\cdot|^{-\gamma}$ is homogeneous of degree $-\gamma$, the change of variables
$y=\sqrt\tau\,y'$ in the convolution integral gives
\begin{equation}\label{Pr.esc}
[S(\tau)|\cdot|^{-\gamma}](x)= \tau^{-\gamma/2}\,F\Big(\frac{x}{\sqrt{\tau}}\Big),
\qquad F:=S(1)|\cdot|^{-\gamma},\qquad \tau>0.
\end{equation}
We first show that $F$ is bounded. Splitting the integral at $|y|=1$ and using
$G_1\le(4\pi)^{-N/2}$ and $\int_{\R^N}G_1=1$, we get, for every $z\in\R^N$,
\[
F(z)\le(4\pi)^{-N/2}\int_{|y|\le1}|y|^{-\gamma}\,dy+\int_{|y|>1}G_1(z-y)\,dy
\le K_1:=\frac{N\omega_N}{N-\gamma}(4\pi)^{-N/2}+1,
\]
which is finite because $\gamma<N$.

Next, we show that $F$ decays like $|z|^{-\gamma}$. Fix $z\neq0$ and split the integral
at $|y|=|z|/2$. If $|y|\le|z|/2$, then $|z-y|\ge|z|/2$, and therefore
$G_1(z-y)\le(4\pi)^{-N/2}e^{-|z|^2/16}$. If $|y|\ge|z|/2$, then
$|y|^{-\gamma}\le2^{\gamma}|z|^{-\gamma}$. It follows that
\[
\begin{aligned}
F(z)&\le(4\pi)^{-N/2}e^{-|z|^2/16}\int_{|y|\le|z|/2}|y|^{-\gamma}\,dy
+2^{\gamma}|z|^{-\gamma}\int_{\R^N}G_1(z-y)\,dy\\
&=\frac{N\omega_N}{N-\gamma}\,2^{\gamma-N}(4\pi)^{-N/2}\,
|z|^{N}e^{-|z|^2/16}\,|z|^{-\gamma}+2^{\gamma}|z|^{-\gamma}.
\end{aligned}
\]
Since $\max_{r>0}r^Ne^{-r^2/16}=(8N/e)^{N/2}$, we deduce that
\[
F(z)\le K_2\,|z|^{-\gamma},\qquad
K_2:=2^{\gamma}+\frac{N\omega_N}{N-\gamma}\,2^{\gamma-N}(4\pi)^{-N/2}
\Big(\frac{8N}{e}\Big)^{N/2}.
\]
We now combine the two bounds. If $|z|\le1$, then $F(z)\le K_1\le 2^\gamma K_1(1+|z|)^{-\gamma}$.
If $|z|\ge1$, then $1+|z|\le2|z|$, and hence $F(z)\le 2^\gamma K_2(1+|z|)^{-\gamma}$.
Setting $K_3:=2^{\gamma}\max\{K_1,K_2\}$, we obtain $F(z)\le K_3(1+|z|)^{-\gamma}$ for all
$z\in\R^N$. Inserting this into \eqref{Pr.esc}, we infer that
\begin{equation}\label{eq:Sweight}
[S(\tau)|\cdot|^{-\gamma}](x)\le K_3\,\tau^{-\gamma/2}
\Big(1+\frac{|x|}{\sqrt\tau}\Big)^{-\gamma}=K_3\,\big(|x|+\sqrt\tau\,\big)^{-\gamma}
\qquad\text{for all }x\in\R^N,\ \tau>0.
\end{equation}

Finally, the substitution $\tau=t-s$ in \eqref{H} and \eqref{eq:Sweight} yield
\[
H(x,t)\le K_3\int_0^t\big(|x|+\sqrt\tau\,\big)^{-\gamma}\,d\tau.
\]
We distinguish two cases. If $|x|\ge\sqrt t$, then $|x|+\sqrt t\le2|x|$, and we get
\[
\int_0^t\big(|x|+\sqrt\tau\,\big)^{-\gamma}\,d\tau\le t\,|x|^{-\gamma}
\le 2^{\gamma}\,t\,\big(|x|+\sqrt t\,\big)^{-\gamma}.
\]
If $|x|<\sqrt t$, then $|x|+\sqrt t<2\sqrt t$, and since $\gamma<2$ we get
\[
\int_0^t\big(|x|+\sqrt\tau\,\big)^{-\gamma}\,d\tau\le\int_0^t\tau^{-\gamma/2}\,d\tau
=\frac{2}{2-\gamma}\,t\,t^{-\gamma/2}
\le\frac{2^{1+\gamma}}{2-\gamma}\,t\,\big(|x|+\sqrt t\,\big)^{-\gamma}.
\]
Since $2/(2-\gamma)\ge1$, both cases give the claim with
$C_1:=2^{1+\gamma}K_3/(2-\gamma)$.
\end{proof}

Since $\big[t(|x|+\sqrt t\,)^{-\gamma}\big]^{1/(1-q)}=W(x,t)$, Lemma~\ref{lemma.H} can be
rewritten as
\begin{equation}\label{eq:HW}
H^{1/(1-q)}\le C_1^{1/(1-q)}\,W\qquad\text{in }\R^N\times(0,\infty).
\end{equation}

The next lemma is the nonlinear companion of Lemma~\ref{lemma.H}. It shows that the
profile $W$ is reproduced, up to a constant, by the nonlinear Duhamel operator, so
that large multiples of $W$ are supersolutions. Recall from \cite[Theorem~1.3(i)]{LM}
that the small multiple $w=c_0W$ is a subsolution.

\begin{lemma}\label{lemma.W}
Let $0<q<1$ and $0<\gamma<\gN$. There exists a constant $C_2=C_2(N,\gamma,q)>0$ such
that
\begin{equation}\label{eq:WW}
\int_0^t\Sg(t-s)\,W^q(s)\,ds\le C_2\,W(t)\qquad\text{in }\R^N,\ \text{for every }t>0.
\end{equation}
Consequently, for every $K\ge C_2^{1/(1-q)}$ the function $KW$ satisfies
$KW(t)\ge\int_0^t\Sg(t-s)(KW)^q(s)\,ds$ for every $t>0$.
\end{lemma}

\begin{proof}
Denote by $I(x,t)$ the left-hand side of \eqref{eq:WW} evaluated at $x$, so that
\[
I(x,t)=\int_0^t\int_{\R^N}G_{t-s}(x-y)\,|y|^{-\gamma}\,s^{q/(1-q)}
\big(|y|+\sqrt s\,\big)^{-\theta q}\,dy\,ds .
\]
We set $m:=\theta q/2$ and $\beta:=q/(1-q)-m=q\alpha/2>0$, and we record the identities
\begin{equation}\label{eq:expo}
\beta+1-\frac\gamma2=\frac\alpha2,\qquad \beta+1+m=\frac1{1-q},\qquad
\gamma+2m=\theta,\qquad \frac1{1-q}-\frac\theta2=\frac\alpha2 .
\end{equation}

We first derive a bound which is uniform in $x$. Using
$(|y|+\sqrt s\,)^{-\theta q}\le s^{-m}$ and then \eqref{eq:Sweight} in the form
$[S(\tau)|\cdot|^{-\gamma}](x)\le K_3\tau^{-\gamma/2}$, we get
\[
I(x,t)\le K_3\int_0^ts^{\beta}(t-s)^{-\gamma/2}\,ds
=K_3\,\mathrm{B}\Big(\beta+1,1-\frac\gamma2\Big)\,t^{\beta+1-\gamma/2}
=K_3\,\mathrm{B}\Big(\beta+1,1-\frac\gamma2\Big)\,t^{\alpha/2},
\]
where $\mathrm{B}$ denotes the Euler beta function, which is finite here because
$\gamma<2$. If $|x|\le\sqrt t$, then $W(x,t)\ge t^{1/(1-q)}(2\sqrt t\,)^{-\theta}
=2^{-\theta}t^{\alpha/2}$ by \eqref{eq:expo}, and therefore
$I(x,t)\le2^{\theta}K_3\,\mathrm{B}(\beta+1,1-\gamma/2)\,W(x,t)$.

Assume now that $|x|>\sqrt t$, and split the integral with respect to $y$ at
$|y|=|x|/2$. If $|y|\ge|x|/2$, then
\[
|y|^{-\gamma}\big(|y|+\sqrt s\,\big)^{-\theta q}\le|y|^{-\gamma-\theta q}=|y|^{-\theta}
\le2^{\theta}|x|^{-\theta},
\]
and since $\int_{\R^N}G_{t-s}=1$, the contribution of this region to $I(x,t)$ is at
most
\[
2^{\theta}|x|^{-\theta}\int_0^ts^{q/(1-q)}\,ds=2^{\theta}(1-q)\,t^{1/(1-q)}|x|^{-\theta}.
\]
If $|y|<|x|/2$, then $|x-y|\ge|x|/2$, hence
$G_{\tau}(x-y)\le(4\pi\tau)^{-N/2}e^{-|x|^2/16\tau}$ with $\tau=t-s$, and we use again
$(|y|+\sqrt s\,)^{-\theta q}\le s^{-m}$. Since
$\int_{|y|<|x|/2}|y|^{-\gamma}dy=\frac{N\omega_N}{N-\gamma}2^{\gamma-N}|x|^{N-\gamma}$,
the contribution of this region is at most
\[
\frac{N\omega_N}{N-\gamma}\,2^{\gamma-N}(4\pi)^{-N/2}\,|x|^{-\gamma}
\int_0^ts^{\beta}\,\Big(\frac{|x|^2}{t-s}\Big)^{N/2}e^{-|x|^2/16(t-s)}\,ds .
\]
Let $\kappa:=\sup_{r>0}r^{N/2+m}e^{-r/16}<\infty$. Then
$r^{N/2}e^{-r/16}\le\kappa\,r^{-m}$ for every $r>0$, and with $r=|x|^2/(t-s)$ we obtain
\[
\Big(\frac{|x|^2}{t-s}\Big)^{N/2}e^{-|x|^2/16(t-s)}\le\kappa\,(t-s)^m|x|^{-2m}
\le\kappa\,t^m|x|^{-2m}.
\]
By \eqref{eq:expo}, the contribution of the region $|y|<|x|/2$ is therefore bounded by
\[
\frac{N\omega_N}{N-\gamma}\,\frac{2^{\gamma-N}(4\pi)^{-N/2}\kappa}{\beta+1}\;
t^{\beta+1+m}\,|x|^{-\gamma-2m}
=\frac{N\omega_N}{N-\gamma}\,\frac{2^{\gamma-N}(4\pi)^{-N/2}\kappa}{\beta+1}\;
t^{1/(1-q)}\,|x|^{-\theta}.
\]
Finally, $|x|>\sqrt t$ implies $|x|+\sqrt t<2|x|$, that is,
$t^{1/(1-q)}|x|^{-\theta}\le2^{\theta}W(x,t)$. Adding the two contributions, we
conclude that \eqref{eq:WW} holds with
\[
C_2:=2^{\theta}\max\Big\{K_3\,\mathrm{B}\Big(\beta+1,1-\frac\gamma2\Big),\
2^{\theta}(1-q)+\frac{N\omega_N}{N-\gamma}\,
\frac{2^{\gamma-N}(4\pi)^{-N/2}\kappa}{\beta+1}\Big\}.
\]
The last assertion follows from the homogeneity of the Duhamel integral, since
$K^qC_2\le K$ if and only if $K\ge C_2^{1/(1-q)}$.
\end{proof}

\subsection{Regularity of the Duhamel integral}
The following property is standard. We include a proof, because it is used in
Section~\ref{sec:zero} to show that a monotone limit of solutions is again a solution
in the sense of Definition~\ref{def:sol}.

\begin{lemma}\label{lemma.reg}
Let $0<\gamma<\gN$ and $T>0$, and let $f:\R^N\times(0,T)\to[0,\infty)$ be a bounded
measurable function. Set $\Psi(0):=0$ and
$\Psi(t):=\int_0^t\Sg(t-\sigma)f(\sigma)\,d\sigma$ for $t\in(0,T]$. Then
$\Psi\in C([0,T];BUC(\R^N))$, and
\begin{equation}\label{eq:regbound}
\|\Psi(t)\|_\infty\le\varkappa\,\|f\|_\infty\,t^{1-\gamma/2}\quad\text{for }t\in[0,T],
\qquad\varkappa:=\frac{2K_3}{2-\gamma},
\end{equation}
where $K_3$ is the constant in \eqref{eq:Sweight}.
\end{lemma}

\begin{proof}
Set $M:=\|f\|_\infty$ and $g(\sigma):=|\cdot|^{-\gamma}f(\sigma)$. By \eqref{eq:Sweight}
we have $0\le S(\tau)g(\sigma)\le MK_3\tau^{-\gamma/2}$ in $\R^N$ for all $\tau>0$ and
$\sigma\in(0,T)$. Hence, for $0\le a<b\le T$,
\begin{equation}\label{eq:Rab}
0\le\int_a^bS(b-\sigma)g(\sigma)\,d\sigma\le MK_3\int_a^b(b-\sigma)^{-\gamma/2}\,d\sigma
=\varkappa\,M\,(b-a)^{1-\gamma/2}.
\end{equation}
Taking $a=0$ gives \eqref{eq:regbound}, and in particular the continuity of $\Psi$ at
$t=0$. Moreover, exactly as in the proof of Lemma~\ref{lemma.inv}(i), Tonelli's theorem
and the semigroup property give, for $0\le s<s+h\le T$,
\begin{equation}\label{eq:Psish}
\Psi(s+h)=S(h)\Psi(s)+\int_s^{s+h}S(s+h-\sigma)g(\sigma)\,d\sigma .
\end{equation}

Let $0<\delta\le s\le T$. Applying \eqref{eq:Psish} with $(s-\delta,\delta)$ in place
of $(s,h)$ and using \eqref{eq:Rab}, we get
\begin{equation}\label{eq:Psidelta}
\big\|\Psi(s)-S(\delta)\Psi(s-\delta)\big\|_\infty\le\varkappa\,M\,\delta^{1-\gamma/2}.
\end{equation}
Since $\Psi(s-\delta)$ is bounded and measurable and $G_\delta$ is smooth with
integrable gradient, the function $S(\delta)\Psi(s-\delta)$ is bounded and Lipschitz
continuous. Letting $\delta\to0$ in \eqref{eq:Psidelta}, we see that $\Psi(s)$ is a
uniform limit of functions of $BUC(\R^N)$, and therefore $\Psi(s)\in BUC(\R^N)$.

It remains to prove the continuity in time on $(0,T]$. Let $0<\delta\le s<s+h\le T$.
By \eqref{eq:Psish}, \eqref{eq:Rab} and \eqref{eq:Psidelta}, and since $S(h)$ is a
contraction on $L^\infty(\R^N)$,
\[
\begin{aligned}
\big\|\Psi(s+h)-\Psi(s)\big\|_\infty
&\le\big\|(S(h)-I)\Psi(s)\big\|_\infty+\varkappa Mh^{1-\gamma/2}\\
&\le\big\|(S(h)-I)S(\delta)\Psi(s-\delta)\big\|_\infty
+2\varkappa M\delta^{1-\gamma/2}+\varkappa Mh^{1-\gamma/2}\\
&\le\|G_{\delta+h}-G_\delta\|_{L^1}\,\varkappa MT^{1-\gamma/2}
+2\varkappa M\delta^{1-\gamma/2}+\varkappa Mh^{1-\gamma/2},
\end{aligned}
\]
where we used $(S(h)-I)S(\delta)\psi=(G_{\delta+h}-G_\delta)\ast\psi$ and
\eqref{eq:regbound}. Given $\varepsilon>0$ and $t_0\in(0,T]$, we first choose
{$\delta\in(0,t_0/2]$} such that $2\varkappa M\delta^{1-\gamma/2}<\varepsilon$, and then
use the continuity of $\tau\mapsto G_\tau$ from $(0,\infty)$ into $L^1(\R^N)$ to find
{$h_0\in(0,t_0/2]$ such that the right-hand side, which does not depend on $s$, is
smaller than $3\varepsilon$ for $0<h<h_0$. Hence
$\|\Psi(s+h)-\Psi(s)\|_\infty<3\varepsilon$ whenever $\delta\le s<s+h\le T$ and
$0<h<h_0$. Since $t_0-h\ge\delta$ for $0<h<h_0$, this applies with $s=t_0$ and with
$s=t_0-h$, and therefore $\Psi$ is continuous at $t_0$.} This completes the proof.
\end{proof}

\section{A comparison lemma}\label{sec:comp}

The next lemma is the key step. It compares a subsolution with a positive
supersolution under the sole assumption that their relative excess is bounded. The
concavity of $s\mapsto s^q$ then improves this bound by the factor $q<1$, and the
excess must vanish.

\begin{lemma}\label{lemma.Comp}
Let $0<q<1$, $\gamma>0$, $T>0$ and $0\le u_0\in C_0(\R^N)$. Let
$u:\R^N\times(0,T]\to[0,\infty)$ and $v:\R^N\times(0,T]\to(0,\infty)$ be measurable
functions satisfying \eqref{eq:subsuper} for every $t\in(0,T]$.
Assume moreover that
\begin{equation}\label{est:Theta}
\Theta:=\sup_{\R^N\times(0,T]} \frac{(u-v)_{+}}{v}< \infty.
\end{equation}
Then $u\leq v$ in $\R^N\times(0,T]$.
\end{lemma}

\begin{proof}
Denote by $J(x,t)$ the Duhamel integral $\int_0^t\Sg(t-s)\,v^q(s)\,ds$ evaluated at $x$.
Since $v$ takes finite values, the second inequality in \eqref{eq:subsuper} shows that
$J$ is finite and that
\begin{equation}\label{eq:J}
J(t)\le v(t)-S(t)u_0\qquad\text{in }\R^N,\ \text{for every }t\in(0,T].
\end{equation}

We next use the concavity of $s\mapsto s^q$. For $a\ge0$ and $b>0$ we have
\begin{equation}\label{eq:tangent}
a^q\le b^q+q\,b^{q-1}(a-b)_+ .
\end{equation}
Indeed, if $a\le b$ this follows from the monotonicity of $s\mapsto s^q$, while if $a>b$
it expresses the fact that the graph of a concave function lies below its tangent line
at $b$. By the definition of $\Theta$, we have $(u-v)_+\le\Theta v$ in
$\R^N\times(0,T]$. Applying \eqref{eq:tangent} with $a=u$ and $b=v$, we obtain
\[
u^q\le v^q+q\,v^{q-1}\,\Theta\, v=(1+q\Theta)\,v^q
\qquad\text{in }\R^N\times(0,T].
\]
Since the Duhamel integral is monotone and positively homogeneous, it follows that
\[
\int_0^t\Sg(t-s)\,u^q(s)\,ds\le(1+q\Theta)\,J(t).
\]
Inserting this into the first
inequality in \eqref{eq:subsuper} and using \eqref{eq:J}, we get
\[
u(t)\le S(t)u_0+(1+q\Theta)\big(v(t)-S(t)u_0\big)
=(1+q\Theta)\,v(t)-q\Theta\, S(t)u_0\le(1+q\Theta)\,v(t),
\]
because $S(t)u_0\ge0$. Hence $(u-v)_+\le q\Theta\,v$ in $\R^N\times(0,T]$. Dividing by
$v>0$ and taking the supremum, we conclude that $\Theta\le q\Theta$. Since $\Theta$ is
finite and $q<1$, it follows that $\Theta=0$, that is, $u\le v$ in $\R^N\times(0,T]$.
\end{proof}

\begin{remark}\rm \label{rem:Theta}
~\begin{enumerate}[label=\rm(\alph*)]
\item The finiteness of $\Theta$ cannot be dropped. Without it, the inequality
$\Theta\le q\Theta$ carries no information. In the proof of Theorem~\ref{thm:super}, $\Theta$ is made
finite by first comparing $u$ with the {regularized} supersolution $Bv+\varepsilon$.
\item Lemma~\ref{lemma.Comp} is related to the comparison principle of Cazenave,
Dickstein and Escobedo \cite[Lemma~6.9 and Theorem~6.2(ii)]{CDE}, proved on a bounded
domain $\Omega$ for $\partial_tu-\Delta u=g(u)$ with Dirichlet boundary conditions.
There, the nonlinearity is required to satisfy the one-sided condition
$g(u)-g(v)\le L_M\,(u-v)/v$ for $0<v\le u\le M$ \cite[(6.6)]{CDE}, which for
$g(s)=s^q$ follows from \eqref{eq:tangent} with $L_M=qM^q$. The supersolution is
bounded from below by a multiple of $\operatorname{dist}(x,\partial\Omega)$, and the
argument is an energy estimate based on Hardy's inequality. Lemma~\ref{lemma.Comp}
works instead with the mild formulation on $\R^N$, and the lower bound on the
supersolution enters through the finiteness of $\Theta$. The mechanism is also the
parabolic counterpart of the classical uniqueness argument for sublinear elliptic
problems and for concave operator equations; see Br\'ezis and Oswald \cite{BO} and
Krasnosel'ski\u{\i} \cite{Kras}. {The same mechanism underlies the ratio argument
of Ioku and Nitta \cite{IokuNitta2026}; see Remark~\ref{rem:IokuNitta}.}
\item Lemma~\ref{lemma.Comp} does not use the specific form of the weight. It holds
for any nonnegative measurable weight in place of $|x|^{-\gamma}$, and for every
$\gamma>0$.
\end{enumerate}
\end{remark}

\section{Proofs of the comparison results}\label{sec:proofs}

\begin{proof}[Proof of Theorem~\ref{thm:super}]
Set $M:=\sup_{\R^N\times(0,T]}u<\infty$. Throughout the proof, all functions are
considered on $\R^N\times(0,T]$, and we write
$J_v(t):=\int_0^t\Sg(t-s)\,v^q(s)\,ds$. Since $v$ takes finite values and satisfies the
second inequality in \eqref{eq:subsuper}, we have
\begin{equation}\label{eq:Jv}
J_v(t)\le v(t)-S(t)u_0,\qquad t\in(0,T].
\end{equation}

\emph{A lower bound for $v$ in terms of $H$.}
By \eqref{eq:HW} and the assumption \eqref{eq:vW}, we have
\begin{equation}\label{eq:HA}
H^{1/(1-q)}\le A\,v\quad\text{in }\R^N\times(0,T],
\qquad A:=\frac{C_1^{1/(1-q)}}{c}.
\end{equation}

\emph{Comparison with a {regularized} supersolution.}
Since $q<1$, we have $B-B^q\to\infty$ as $B\to\infty$. We may therefore fix $B\ge1$
such that
\begin{equation}\label{eq:B}
B^q+(1-q)A\le B.
\end{equation}
For $\varepsilon>0$, we set $v_\varepsilon:=Bv+\varepsilon$, and we claim that
$v_\varepsilon$ satisfies the second inequality in \eqref{eq:subsuper}. The
subadditivity of $s\mapsto s^q$ on $[0,\infty)$ gives
$v_\varepsilon^q\le B^qv^q+\varepsilon^q$. Since $\Sg(t-s)\varepsilon^q=\varepsilon^q
S(t-s)|\cdot|^{-\gamma}$, it follows from \eqref{H} and \eqref{eq:Jv} that
\[
\int_0^t\Sg(t-s)\,v_\varepsilon^q(s)\,ds\le B^qJ_v(t)+\varepsilon^qH(t)
\le B^q\big(v(t)-S(t)u_0\big)+\varepsilon^qH(t).
\]
By Young's inequality $ab\le q\,a^{1/q}+(1-q)\,b^{1/(1-q)}$, valid for $a,b\ge0$, and by
\eqref{eq:HA}, we have
\[
\varepsilon^qH\le q\,\varepsilon+(1-q)\,H^{1/(1-q)}\le q\,\varepsilon+(1-q)A\,v.
\]
Moreover, $(1-B^q)S(t)u_0\le0$, because $B\ge1$ and the semigroup is positivity
preserving. Combining these estimates with \eqref{eq:B}, we obtain
\[
\begin{aligned}
S(t)u_0+\int_0^t\Sg(t-s)\,v_\varepsilon^q(s)\,ds
&\le(1-B^q)S(t)u_0+\big(B^q+(1-q)A\big)v(t)+q\,\varepsilon\\
&\le B\,v(t)+\varepsilon=v_\varepsilon(t).
\end{aligned}
\]
This proves the claim. In addition, $v_\varepsilon\ge\varepsilon$ and $0\le u\le M$, so that
\[
\sup_{\R^N\times(0,T]}\frac{(u-v_\varepsilon)_+}{v_\varepsilon}\le\frac{M}{\varepsilon}<\infty.
\]
Lemma~\ref{lemma.Comp}, applied to $u$ and $v_\varepsilon$, yields $u\le Bv+\varepsilon$ in
$\R^N\times(0,T]$. Letting $\varepsilon\to0$, we conclude that
\begin{equation}\label{eq:uBv}
u\le B\,v\qquad\text{in }\R^N\times(0,T].
\end{equation}

\emph{Conclusion.}
The function $v$ is positive and satisfies the second inequality in
\eqref{eq:subsuper}, and by \eqref{eq:uBv} we have $(u-v)_+\le(B-1)v$, so that the
quantity $\Theta$ in \eqref{est:Theta} is at most $B-1<\infty$.
Lemma~\ref{lemma.Comp} therefore gives $u\le v$ in $\R^N\times(0,T]$.
\end{proof}

\begin{proof}[Proof of Theorem~\ref{thm:comparison}]
Fix $T>0$. Since $v$ is a solution with initial data $v_0\ge u_0$ and the semigroup is
positivity preserving, we have
\[
v(t)=S(t)v_0+\int_0^t\Sg(t-s)\,v^q(s)\,ds\ge S(t)u_0+\int_0^t\Sg(t-s)\,v^q(s)\,ds,
\]
so that $v$ satisfies the second inequality in \eqref{eq:subsuper} with initial datum
$u_0$. Since $v_0\not\equiv0$, Proposition~\ref{prop:lower} shows that
$v\ge w=c_0W>0$ in $\R^N\times(0,T]$, which is \eqref{eq:vW} with $c=c_0$. Finally, $u$
is bounded on $\R^N\times(0,T]$ and satisfies the first inequality in
\eqref{eq:subsuper}. Theorem~\ref{thm:super} gives $u\le v$ in $\R^N\times(0,T]$, and
since $T>0$ is arbitrary, the proof is complete. Note that in this case the constant
$A$ in \eqref{eq:HA} equals $[C_1/((1-q)\eta_0)]^{1/(1-q)}$ and depends only on $N$,
$\gamma$ and $q$.
\end{proof}

\begin{proof}[Proof of Corollary~\ref{cor:uniq}]
Let $u$ be a nonnegative solution of \eqref{eq:main} with initial data $u_0$. Since
$u\in C([0,\infty);C_0(\R^N))$, the function $(x,t)\mapsto u(x,t)$ is continuous on
$\R^N\times[0,\infty)$, hence measurable. Moreover, $t\mapsto\|u(t)\|_\infty$ is
continuous, and therefore $u$ is bounded on $\R^N\times[0,T]$ for every $T>0$. Finally,
$u$ satisfies \eqref{eq:subsol} with equality. Thus $u$ is admissible in
Theorem~\ref{thm:comparison}.

Assertion (ii) follows directly from Theorem~\ref{thm:comparison} on
$\R^N\times(0,\infty)$, while at $t=0$ it reduces to $u_0\le v_0$. For assertion (i),
existence is provided by \cite[Theorem~1.1]{LM}. If $u$ and $v$ are two nonnegative
solutions with the same initial data $u_0\not\equiv0$, then (ii) applied with
$v_0=u_0$ gives $u\le v$. Exchanging the roles of $u$ and $v$, we also get $v\le u$, and
we conclude that $u=v$.
\end{proof}

\section{Solutions emanating from zero}\label{sec:zero}

This section is devoted to the proof of Theorem~\ref{thm:zero}. The strategy follows
the scheme of \cite[Corollary~2.6]{AE}: a solution which leaves zero at time $\tau$
satisfies the lower bound \eqref{eq:lower} from that time on, and the lower bound
restores uniqueness. The explicit solution available for $\gamma=0$ is replaced here
by a solution constructed as a monotone limit, and the explicit upper bound available
in \cite{AE} for solutions with zero initial data is replaced by Theorem~\ref{thm:super}.

\begin{proof}[Proof of Theorem~\ref{thm:zero}]
\emph{Construction of $U$.}
Let $\varphi(x):=e^{-|x|^2}$, and for $n\ge1$ let $u_n$ be the unique nonnegative
solution of \eqref{eq:main} with initial data $\varphi/n$, given by
Corollary~\ref{cor:uniq}(i). By Corollary~\ref{cor:uniq}(ii), the sequence $(u_n)$ is
nonincreasing, and by Proposition~\ref{prop:lower} we have $u_n\ge c_0W$ for every $n$.
We define
\[
U(x,t):=\lim_{n\to\infty}u_n(x,t)=\inf_{n\ge1}u_n(x,t),\qquad (x,t)\in\R^N\times[0,\infty).
\]
Then $U$ is measurable, $U(0)=0$ and
\begin{equation}\label{eq:Usandwich}
c_0\,W\le U\le u_1\qquad\text{in }\R^N\times[0,\infty).
\end{equation}
Fix $t>0$ and $x\in\R^N$. We have $0\le[S(t)\varphi](x)/n\to0$, and $u_n^q\downarrow U^q$
with $u_n^q\le u_1^q$. Since the Duhamel integral of $u_1^q$ at $(x,t)$ is finite, being
bounded by $u_1(x,t)$, the dominated convergence theorem allows us to pass to the
limit in the integral equation satisfied by $u_n$, and we obtain
\begin{equation}\label{eq:Uduhamel}
U(t)=\int_0^t\Sg(t-\sigma)\,U^q(\sigma)\,d\sigma\qquad\text{in }\R^N,\ \text{for every }t\ge0.
\end{equation}
Let $T>0$. By \eqref{eq:Usandwich}, the function $U^q$ is bounded and measurable on
$\R^N\times(0,T)$, and Lemma~\ref{lemma.reg} together with \eqref{eq:Uduhamel} shows
that $U\in C([0,T];BUC(\R^N))$. Moreover, $0\le U(t)\le u_1(t)$ and $u_1(t)\in C_0(\R^N)$,
so that $U(t)$ vanishes at infinity. Hence $U\in C([0,\infty);C_0(\R^N))$, and $U$ is a
nonnegative solution of \eqref{eq:main} with $u_0\equiv0$ which is positive in
$\R^N\times(0,\infty)$ and satisfies the first inequality in \eqref{eq:Ubounds}.

\emph{Proof of (iii).}
By \eqref{eq:Uduhamel} and \eqref{eq:Usandwich}, the function $U$ is positive in
$\R^N\times(0,\infty)$, satisfies the second inequality in \eqref{eq:subsuper} with
$u_0\equiv0$, and satisfies \eqref{eq:vW} with $c=c_0$. Let $u$ be a nonnegative
measurable function which is bounded on $\R^N\times(0,T]$ for every $T>0$ and satisfies
\eqref{eq:subsol} with $u_0\equiv0$. Theorem~\ref{thm:super}, applied with $v=U$ on
each interval $(0,T]$, gives $u\le U$. Since $U$ itself belongs to this class of
functions, it is its largest element. If $v$ is a nonnegative solution with initial
data $v_0\not\equiv0$, then Theorem~\ref{thm:comparison}, applied with $u=U$ and
$u_0\equiv0$, gives $U\le v$.

\emph{Proof of (ii).}
We first check that $U_\tau$ is a solution with zero initial data for every
$\tau\ge0$. Clearly $U_\tau\in C([0,\infty);C_0(\R^N))$ and $U_\tau(t)=0$ for
$t\in[0,\tau]$, so that \eqref{eq:duhamel} holds trivially on $[0,\tau]$. For
$t>\tau$, since $U_\tau^q(\sigma)=0$ for $\sigma\le\tau$, the change of variable
$\sigma=\tau+\sigma'$ and \eqref{eq:Uduhamel} give
\[
\int_0^t\Sg(t-\sigma)\,U_\tau^q(\sigma)\,d\sigma
=\int_0^{t-\tau}\Sg(t-\tau-\sigma')\,U^q(\sigma')\,d\sigma'=U(t-\tau)=U_\tau(t).
\]

Conversely, let $u$ be a nonnegative solution of \eqref{eq:main} with $u_0\equiv0$ and
$u\not\equiv0$, and consider the set
\[
Z:=\{t\ge0:\ u(t)\equiv0\}.
\]
Then $0\in Z$, and $Z$ is closed because $u\in C([0,\infty);C_0(\R^N))$. Let
$s\ge0$ with $s\notin Z$. By Lemma~\ref{lemma.inv}(i), the function $u(\cdot+s)$ is a
nonnegative solution of \eqref{eq:main} with initial data $u(s)\not\equiv0$, and
Proposition~\ref{prop:lower} gives
\begin{equation}\label{eq:shiftlower}
u(x,s+t)\ge c_0\,W(x,t)\qquad\text{for all }x\in\R^N,\ t\ge0.
\end{equation}
In particular, $u(s+t)\not\equiv0$ for every $t>0$, that is, $(s,\infty)\cap Z=\emptyset$.
Consequently, $Z$ is a closed interval containing $0$, and since $u\not\equiv0$ we
have $Z=[0,\tau]$ for some $\tau\in[0,\infty)$. Now \eqref{eq:shiftlower} holds for
every $s>\tau$. Letting $s\downarrow\tau$ and using the continuity of $u$, we obtain
\begin{equation}\label{eq:tildelower}
\tilde u(x,t):=u(x,\tau+t)\ge c_0\,W(x,t)\qquad\text{for all }x\in\R^N,\ t\ge0.
\end{equation}
By Lemma~\ref{lemma.inv}(i) again, $\tilde u$ is a nonnegative solution of
\eqref{eq:main} with initial data $u(\tau)\equiv0$. In particular, $\tilde u$ is
bounded on $\R^N\times(0,T]$ for every $T>0$ and satisfies \eqref{eq:subsol} with
$u_0\equiv0$, so that $\tilde u\le U$ by (iii). On the other hand, by
\eqref{eq:tildelower} the function $\tilde u$ is positive in $\R^N\times(0,\infty)$ and
satisfies \eqref{eq:vW} with $c=c_0$, and it satisfies the second inequality in
\eqref{eq:subsuper} with $u_0\equiv0$. Theorem~\ref{thm:super}, applied with the
subsolution $U$ and the supersolution $\tilde u$ on each interval $(0,T]$, gives
$U\le\tilde u$. Therefore $\tilde u=U$, which means that $u(t)=U(t-\tau)$ for
$t\ge\tau$. Since $u(t)\equiv0$ for $t\in[0,\tau]$, we conclude that $u=U_\tau$.

\emph{Proof of (i).}
We have already shown that $U$ is a nonnegative solution with zero initial data such
that $U(t)\not\equiv0$ for $t>0$. If $u$ is another solution with these properties,
then $u=U_\tau$ for some $\tau\ge0$ by (ii), and $\tau=0$ because $U_\tau(t)\equiv0$
for $t\le\tau$. This proves the uniqueness assertion.

Let $\lambda>0$. By Lemma~\ref{lemma.inv}(ii), the function
$(x,t)\mapsto\lambda^{-\alpha}U(\lambda x,\lambda^2t)$ is a nonnegative solution with
zero initial data which does not vanish identically at any positive time. By
uniqueness,
\[
\lambda^{-\alpha}\,U(\lambda x,\lambda^2t)=U(x,t)\qquad\text{for all }x\in\R^N,\ t\ge0,\
\lambda>0 .
\]
Replacing $x$ by $x/\sqrt t$ and $t$ by $1$, and choosing $\lambda=\sqrt t$, we obtain
\eqref{eq:selfsim}. In the same way, Lemma~\ref{lemma.inv}(iii) and uniqueness give
$U(Rx,t)=U(x,t)$ for every orthogonal matrix $R$, so that $U(\cdot,t)$ is radially
symmetric.

It remains to prove the upper bound in \eqref{eq:Ubounds}. Set
$K_0:=C_2^{1/(1-q)}$, where $C_2$ is the constant of Lemma~\ref{lemma.W}. By that
lemma, the function $v:=K_0W$ is positive in $\R^N\times(0,\infty)$ and satisfies the
second inequality in \eqref{eq:subsuper} with $u_0\equiv0$, and it obviously satisfies
\eqref{eq:vW} with $c=K_0$. Since $U$ is bounded on $\R^N\times(0,T]$ for every $T>0$
and satisfies \eqref{eq:subsol} with $u_0\equiv0$, Theorem~\ref{thm:super} gives
$U\le K_0W$. The bounds on $\Phi$ follow by taking $t=1$ in \eqref{eq:Ubounds}, since
$W(z,1)=(1+|z|)^{-\theta}$.

\emph{Proof of (iv).}
By Corollary~\ref{cor:uniq}, the solutions $v_n$ are well defined and
$v_n\ge v_{n+1}$, and by (iii) we have $v_n\ge U$ for every $n$. Let
$V:=\lim_{n\to\infty}v_n\ge U$. Since $S(t)v_{0,n}\to0$ pointwise by dominated
convergence, the argument used in the construction of $U$ shows that $V$ is
measurable, bounded on $\R^N\times(0,T]$ for every $T>0$, and satisfies
\eqref{eq:subsol} with $u_0\equiv0$, with equality. Hence $V\le U$ by (iii), and
therefore $V=U$.

To prove that the convergence is uniform, fix $T>0$ and $\varepsilon>0$. The set
$\{v_1(t):t\in[0,T]\}$ is compact in $C_0(\R^N)$, as the continuous image of a compact
interval. Covering it by finitely many balls of radius $\varepsilon/2$, we find $R>0$
such that $v_1(x,t)<\varepsilon$ for all $|x|\ge R$ and $t\in[0,T]$, and hence
$0\le v_n-U\le v_1<\varepsilon$ in this region for every $n$. On the compact set
$\{|x|\le R\}\times[0,T]$, the continuous functions $v_n-U$ decrease pointwise to
zero, and Dini's theorem shows that the convergence is uniform. Therefore
$\sup_{\R^N\times[0,T]}(v_n-U)\to0$ as $n\to\infty$.
\end{proof}

\begin{remark}\rm \label{rem:zero}
~\begin{enumerate}[label=\rm(\alph*)]
\item The proof of (ii) shows slightly more than stated: a nonnegative solution with
an arbitrary initial data $u_0\in C_0(\R^N)$ either vanishes identically on
$\R^N\times[0,\infty)$, which forces $u_0\equiv0$, or is positive on
$\R^N\times(\tau,\infty)$ for some $\tau\ge0$, with $\tau=0$ if $u_0\not\equiv0$.
\item The existence scheme of \cite{LM} approximates the initial data from above.
{For $u_0\equiv0$ we therefore expect the solution constructed in
\cite[Theorem~1.1]{LM} to coincide with the maximal solution $U$; we have not
verified this point, which is not needed here.}
\item Theorem~\ref{thm:zero}(iv) expresses a form of continuous dependence at the
degenerate data $u_0\equiv0$: among the infinitely many solutions leaving zero, the
solutions with small nontrivial data select the maximal one. {For $\gamma=0$ and
constant data $\varepsilon>0$, the same selection can be seen explicitly: the solution is
$[\varepsilon^{1-q}+(1-q)t]^{1/(1-q)}$, which converges to $[(1-q)t]^{1/(1-q)}$ as
$\varepsilon\to0$.}
\end{enumerate}
\end{remark}

\section{Bounds and large-time behavior}\label{sec:asym}

\begin{proof}[Proof of Theorem~\ref{thm:asym}]
\emph{Proof of (i).}
The inequality $U\le u$ is contained in Theorem~\ref{thm:zero}(iii). For the upper
bound, set $a:=\|u_0\|_\infty>0$. Since $K-K^qC_2\to\infty$ as $K\to\infty$, we may fix
$K_1=K_1(N,\gamma,q)>0$ such that
\begin{equation}\label{eq:K1}
K_1^q\,C_2+(1-q)\big(2^qC_1\big)^{1/(1-q)}\le K_1,
\end{equation}
where $C_1$ and $C_2$ are the constants of Lemmas~\ref{lemma.H} and~\ref{lemma.W}. We
claim that $v:=2a+K_1W$ satisfies the second inequality in \eqref{eq:subsuper} for
every $t>0$. By subadditivity, $v^q\le(2a)^q+K_1^qW^q$, and therefore, by \eqref{H} and
Lemma~\ref{lemma.W},
\[
\int_0^t\Sg(t-s)\,v^q(s)\,ds\le(2a)^qH(t)+K_1^q\,C_2\,W(t).
\]
Young's inequality and \eqref{eq:HW} give
\[
(2a)^qH=a^q\,\big(2^qH\big)\le q\,a+(1-q)\big(2^qH\big)^{1/(1-q)}
\le q\,a+(1-q)\big(2^qC_1\big)^{1/(1-q)}\,W .
\]
Since $S(t)u_0\le a$, we deduce from \eqref{eq:K1} that
\[
S(t)u_0+\int_0^t\Sg(t-s)\,v^q(s)\,ds\le a+q\,a+K_1\,W(t)\le2a+K_1W(t)=v(t),
\]
which proves the claim. The function $v$ is positive and satisfies \eqref{eq:vW} with
$c=K_1$, while $u$ is bounded on $\R^N\times(0,T]$ for every $T>0$ and satisfies the
first inequality in \eqref{eq:subsuper}. Theorem~\ref{thm:super} gives $u\le v$ on
$\R^N\times(0,\infty)$, and the inequality is trivial at $t=0$. This proves
\eqref{eq:twosided}. The bounds on $\|u(t)\|_\infty$ follow from
$u(0,t)\ge U(0,t)\ge c_0W(0,t)=c_0t^{\alpha/2}$ and from $\sup_xW(x,t)=t^{\alpha/2}$.

\emph{Proof of (ii).}
We introduce the radial majorant of the initial datum,
\[
\psi(x):=\sup_{|y|\ge|x|}u_0(y),\qquad x\in\R^N .
\]
The function $\psi$ is radially nonincreasing, $u_0\le\psi$, $\psi(0)=\|u_0\|_\infty>0$
and $\psi(x)\to0$ as $|x|\to\infty$. Moreover, $\psi$ is continuous. Indeed, let
$\rho$ be a modulus of continuity of the uniformly continuous function $u_0$, and let
$x,x'\in\R^N$. For every $y$ with $|y|\ge|x|$ there exists $y'$ with $|y'|\ge|x'|$ and
$|y-y'|\le\big||x|-|x'|\big|$, so that $u_0(y)\le\psi(x')+\rho(|x-x'|)$; taking the
supremum over $y$ and exchanging the roles of $x$ and $x'$ gives
$|\psi(x)-\psi(x')|\le\rho(|x-x'|)$. Thus $\psi\in C_0(\R^N)$.

For $\lambda\ge1$, let $u_\lambda(x,t):=\lambda^{-\alpha}u(\lambda x,\lambda^2t)$. By
Lemma~\ref{lemma.inv}(ii), $u_\lambda$ is the solution of \eqref{eq:main} with initial
data $\lambda^{-\alpha}u_0(\lambda\,\cdot)$, and since $\psi$ is radially nonincreasing
and $\lambda\ge1$,
\[
0\le\lambda^{-\alpha}u_0(\lambda x)\le\lambda^{-\alpha}\psi(\lambda x)\le\lambda^{-\alpha}\psi(x),
\qquad x\in\R^N .
\]
Let $V_\lambda$ denote the solution of \eqref{eq:main} with initial datum
$\lambda^{-\alpha}\psi\not\equiv0$. Corollary~\ref{cor:uniq}(ii) and
Theorem~\ref{thm:zero}(iii) give
\begin{equation}\label{eq:squeeze}
U\le u_\lambda\le V_\lambda\qquad\text{in }\R^N\times[0,\infty),
\end{equation}
and $V_\lambda\le V_\mu$ for $\lambda\ge\mu\ge1$. By Theorem~\ref{thm:zero}(iv),
applied to the sequence of initial data $n^{-\alpha}\psi$, we have
$\sup_{\R^N\times[0,1]}(V_n-U)\to0$ as $n\to\infty$, and since
$\lambda\mapsto\sup_{\R^N\times[0,1]}(V_\lambda-U)$ is nonincreasing, the same holds as
$\lambda\to\infty$ along the real numbers. Evaluating \eqref{eq:squeeze} at $t=1$, we
obtain
\[
\lim_{\lambda\to\infty}\ \sup_{x\in\R^N}\big|\lambda^{-\alpha}u(\lambda x,\lambda^2)-U(x,1)\big|=0 .
\]
By the self-similarity of $U$, we have $U(x,1)=\lambda^{-\alpha}U(\lambda x,\lambda^2)$.
Setting $t=\lambda^2$ and $y=\lambda x$, we conclude that
\[
t^{-\alpha/2}\,\sup_{y\in\R^N}\big|u(y,t)-U(y,t)\big|\longrightarrow0
\qquad\text{as }t\to\infty,
\]
which is \eqref{eq:asym}. The equivalent formulation follows from \eqref{eq:selfsim}.
\end{proof}

\begin{remark}\rm \label{rem:asym}
The conclusion of Theorem~\ref{thm:asym}(ii) also holds for the solutions $U_\tau$,
$\tau>0$, emanating from zero. Indeed, by Lemma~\ref{lemma.inv}(i) the function
$U(\cdot+\tau)$ is the solution with initial data $U(\tau)\not\equiv0$, so that
$t^{-\alpha/2}\|U(t+\tau)-U(t)\|_\infty\to0$ as $t\to\infty$ by \eqref{eq:asym}.
Replacing $t$ by $t-\tau$ and using $(t-\tau)/t\to1$, we get
$t^{-\alpha/2}\|U(t)-U_\tau(t)\|_\infty\to0$. Hence every nontrivial nonnegative
solution of \eqref{eq:main}, whatever its initial datum, behaves like $U$ as
$t\to\infty$.
\end{remark}

\section{Concluding remarks}\label{sec:concl}

We close with a few comments on the scope of the method.

The comparison argument is not tied to the power weight. Let $a\ge0$ be a measurable
weight, replace $|x|^{-\gamma}$ by $a$ in \eqref{eq:main}, and set
$H_a(t)=\int_0^tS(\tau)a\,d\tau$. Suppose that $H_a$ is finite on
$\R^N\times(0,T]$ and that a positive supersolution $v$ satisfies
$H_a^{1/(1-q)}\le A\,v$ on $\R^N\times(0,T]$ for some constant $A$. Then the proof of
Theorem~\ref{thm:super} applies verbatim, since Lemma~\ref{lemma.Comp} holds for any
nonnegative weight. In this sense, the problem of uniqueness for concave reactions
reduces to a lower bound of Aguirre--Escobedo type, expressed through the heat flow
of the weight. {The description of the solutions emanating from zero given in
Theorem~\ref{thm:zero}(ii) should extend to such weights as soon as every nontrivial
solution satisfies a lower bound of this type by a fixed function $\underline w$ with
$H_a^{1/(1-q)}\le A\,\underline w$; the self-similar structure and the bounds
\eqref{eq:Ubounds} are of course specific to homogeneous weights.}

We believe that the localized weights
$a=\mathbf{1}_{B_L}|x|^{-\gamma}$, with $L<\infty$, considered by
Ioku and Nitta \cite{IokuNitta2026}, fit within our framework, provided
that every nontrivial solution satisfies $H_a^{1/(1-q)}\le A\,v$ throughout the whole space. In particular, this estimate must remain
valid in the exterior region $\{|x|>L\}$, where the reaction term vanishes. The lower bound of
\cite[Lemma~2.2]{IokuNitta2026} is stated only in $B_L$. Comparing the Gaussian tails
of $H_a$ and of the Duhamel term of $v$ suggests that the required bound also holds
outside $B_L$, which would remove the condition $\gamma<N(1-q)$ for localized weights as
well. We have not checked the details.

Boundedness of the subsolution $u$ enters only through the estimate
$(u-v_\varepsilon)_+/v_\varepsilon\le M/\varepsilon$. It would be interesting to know whether
uniqueness persists for unbounded initial data, for instance in the uniformly local
Lebesgue spaces considered in \cite{CCL}, where existence is available, or in weighted
spaces of the type used in \cite{AE}. In such classes one expects, as in \cite{AE}, a
richer family of self-similar solutions and a large-time behavior depending on the
behavior of the initial data at infinity.

The profile $\Phi$ of the self-similar solution $U$ deserves a separate study:
monotonicity, regularity at the origin, the precise behavior at infinity, and the
uniqueness of positive solutions of the elliptic equation in
Remark~\ref{rem:main}(e) are open. It would also be natural to investigate whether
the trichotomy of Theorem~\ref{thm:zero}(ii) persists for weakly coupled systems with
singular weights, in the spirit of \cite{EH93,DE}.

Finally, the order structure is essential in our argument. The uniqueness of
sign-changing solutions for the nonlinearity $|x|^{-\gamma}|u|^{q-1}u$, as well as the
borderline case $\gamma=2$ with $N\ge3$, remain open.

\section*{Declarations}
\noindent\textbf{Funding.} This research received no external funding.

\noindent\textbf{Conflict of interest.} The authors declare that they have no conflicts of interest.

\noindent\textbf{Data availability.} No datasets were generated or analyzed during the current study.


\begin{thebibliography}{99}

\bibitem{AgarwalLakshmikantham} R.~P.~Agarwal and V.~Lakshmikantham,
\emph{Uniqueness and Nonuniqueness Criteria for Ordinary Differential Equations},
Ser.\ Real Anal.\ \textbf{6}, World Scientific, Singapore, 1993.

\bibitem{AE} J.~Aguirre and M.~Escobedo,
\emph{A Cauchy problem for $u_t-\Delta u=u^p$ with $0<p<1$. Asymptotic behaviour of
solutions}, Ann.\ Fac.\ Sci.\ Toulouse Math.\ (5) \textbf{8} (1986/87), no.~2, 175--203.

\bibitem{BG} P.~Baras and J.~A.~Goldstein,
\emph{The heat equation with a singular potential},
Trans.\ Amer.\ Math.\ Soc.\ \textbf{284} (1984), no.~1, 121--139.

\bibitem{BSTW} B.~Ben Slimene, S.~Tayachi and F.~B.~Weissler,
\emph{Well-posedness, global existence and large time behavior for Hardy--H\'enon
parabolic equations}, Nonlinear Anal.\ \textbf{152} (2017), 116--148.

\bibitem{BC} H.~Br\'ezis and T.~Cazenave,
\emph{A nonlinear heat equation with singular initial data},
J.\ Anal.\ Math.\ \textbf{68} (1996), 277--304.

\bibitem{BO} H.~Br\'ezis and L.~Oswald,
\emph{Remarks on sublinear elliptic equations},
Nonlinear Anal.\ \textbf{10} (1986), no.~1, 55--64.

\bibitem{CM} X.~Cabr\'e and Y.~Martel,
\emph{Existence versus explosion instantan\'ee pour des \'equations de la chaleur
lin\'eaires avec potentiel singulier}, C.\ R.\ Acad.\ Sci.\ Paris S\'er.\ I Math.\
\textbf{329} (1999), no.~11, 973--978.

\bibitem{CCL} B.~Carhuas-Torre, R.~Castillo and M.~Loayza,
\emph{The Hardy parabolic problem with initial data in uniformly local Lebesgue spaces},
J.\ Differential Equations \textbf{424} (2025), 438--462.

\bibitem{CGL} R.~Castillo, O.~Guzm\'an-Rea and M.~Loayza,
\emph{On the local existence for Hardy parabolic equations with singular initial data},
J.\ Math.\ Anal.\ Appl.\ \textbf{510} (2022), no.~2, Paper No.~126022.

\bibitem{CDE} T.~Cazenave, F.~Dickstein and M.~Escobedo,
\emph{A semilinear heat equation with concave-convex nonlinearity},
Rend.\ Mat.\ Appl.\ (7) \textbf{19} (1999), no.~2, 211--242.

\bibitem{CITT} N.~Chikami, M.~Ikeda, K.~Taniguchi and S.~Tayachi,
\emph{Unconditional uniqueness and non-uniqueness for Hardy--H\'enon parabolic
equations}, Math.\ Ann.\ \textbf{390} (2024), no.~3, 3765--3825.

\bibitem{CRB} J.~W.~Cholewa and A.~Rodriguez-Bernal,
\emph{Self-similarity in homogeneous stationary and evolution problems},
J.\ Evol.\ Equ.\ \textbf{23} (2023), no.~2, Paper No.~42.

\bibitem{DE} F.~Dickstein and M.~Escobedo,
\emph{A maximum principle for semilinear parabolic systems and applications},
Nonlinear Anal.\ \textbf{45} (2001), 825--837.

\bibitem{EH} M.~Escobedo and M.~A.~Herrero,
\emph{A uniqueness result for a semilinear reaction-diffusion system},
Proc.\ Amer.\ Math.\ Soc.\ \textbf{112} (1991), no.~1, 175--185.

\bibitem{EH93} M.~Escobedo and M.~A.~Herrero,
\emph{A semilinear parabolic system in a bounded domain},
Ann.\ Mat.\ Pura Appl.\ (4) \textbf{165} (1993), 315--336.

\bibitem{FW} H.~Fujita and S.~Watanabe,
\emph{On the uniqueness and non-uniqueness of solutions of initial value problems for
some quasi-linear parabolic equations}, Comm.\ Pure Appl.\ Math.\ \textbf{21} (1968),
631--652.

\bibitem{HW} A.~Haraux and F.~B.~Weissler,
\emph{Nonuniqueness for a semilinear initial value problem},
Indiana Univ.\ Math.\ J.\ \textbf{31} (1982), 167--189.

\bibitem{Henry} D.~Henry,
\emph{Geometric Theory of Semilinear Parabolic Equations},
Lecture Notes in Math.\ \textbf{840}, Springer, Berlin, 1981.

\bibitem{HiTs} D.~Hirata and M.~Tsutsumi,
\emph{On the well-posedness of a linear heat equation with a critical singular
potential}, Differential Integral Equations \textbf{14} (2001), no.~1, 1--18.

\bibitem{HS} K.~Hisa and M.~Sier\.z\c{e}ga,
\emph{Existence and nonexistence of solutions to the Hardy parabolic equation},
Funkcial.\ Ekvac.\ \textbf{67} (2024), no.~2, 149--174.

\bibitem{HT} K.~Hisa and J.~Takahashi,
\emph{Optimal singularities of initial data for solvability of the Hardy parabolic
equation}, J.\ Differential Equations \textbf{296} (2021), 822--848.

\bibitem{IokuNitta2026} N.~Ioku and Y.~Nitta, \emph{Uniqueness of global solutions for the
semilinear heat equation with localized concave nonlinearities}, Bull.\ Lond.\ Math.\ Soc.\
\textbf{58} (2026), {no.~12, Paper No.~e70552.}

\bibitem{Kras} M.~A.~Krasnosel'ski\u{\i},
\emph{Positive Solutions of Operator Equations}, P.~Noordhoff, Groningen, 1964.

\bibitem{Loayza06} M.~Loayza,
\emph{The heat equation with singular nonlinearity and singular initial data},
J.\ Differential Equations \textbf{229} (2006), no.~2, 509--528.

\bibitem{LM} M.~Loayza and M.~Majdoub,
\emph{An inhomogeneous nonlinear parabolic equation with concave nonlinearity},
Tunis.\ J.\ Math.\ \textbf{8} (2026), no.~1, 41--58.

\bibitem{Maj21} M.~Majdoub,
\emph{Well-posedness and blow-up for an inhomogeneous semilinear parabolic equation},
Differ.\ Equ.\ Appl.\ \textbf{13} (2021), no.~1, 85--100.

\bibitem{Maj23} M.~Majdoub,
\emph{On the Fujita exponent for a Hardy--H\'enon equation with a spatial-temporal
forcing term}, La Matematica \textbf{2} (2023), no.~2, 340--361.

\bibitem{MM21} M.~Majdoub and E.~Mliki,
\emph{Well-posedness for Hardy--H\'enon parabolic equations with fractional Brownian
noise}, Anal.\ Math.\ Phys.\ \textbf{11} (2021), no.~1, Paper No.~20.

\bibitem{Nagumo} M.~Nagumo,
\emph{Eine hinreichende Bedingung f\"ur die Unit\"at der L\"osung von
Differentialgleichungen erster Ordnung},
Jpn.\ J.\ Math.\ \textbf{3} (1926), 107--112.

\bibitem{Osgood} W.~F.~Osgood,
\emph{Beweis der Existenz einer L\"osung der Differentialgleichung
$\frac{dy}{dx}=f(x,y)$ ohne Hinzunahme der Cauchy--Lipschitz'schen Bedingung},
Monatsh.\ Math.\ Phys.\ \textbf{9} (1898), 331--345.

\bibitem{QS} P.~Quittner and P.~Souplet,
\emph{Superlinear Parabolic Problems. Blow-up, Global Existence and Steady States},
2nd ed., Birkh\"auser Adv.\ Texts, Birkh\"auser/Springer, Cham, 2019.

\bibitem{Tayachi} S.~Tayachi,
\emph{Uniqueness and non-uniqueness of solutions for critical Hardy--H\'enon parabolic
equations}, J.\ Math.\ Anal.\ Appl.\ \textbf{488} (2020), no.~1, Paper No.~123976.

\bibitem{VZ} J.~L.~V\'azquez and E.~Zuazua,
\emph{The Hardy inequality and the asymptotic behaviour of the heat equation with an
inverse-square potential}, J.\ Funct.\ Anal.\ \textbf{173} (2000), no.~1, 103--153.

\bibitem{W80} F.~B.~Weissler,
\emph{Local existence and nonexistence for semilinear parabolic equations in $L^p$},
Indiana Univ.\ Math.\ J.\ \textbf{29} (1980), no.~1, 79--102.

\end{thebibliography}
\end{document}